\documentclass[11pt]{article}
\usepackage{amsmath,amssymb,amsthm}
\usepackage{authblk}
\usepackage{xcolor}
\usepackage[a4paper,top=2.5cm,bottom=2.5cm,left=3cm,right=2.5cm]{geometry}

\newtheorem{theorem}{Theorem}
\newtheorem{lemma}[theorem]{Lemma}

\theoremstyle{remark}
\newtheorem{remark}[theorem]{Remark}

\newcommand{\KK}{\mathcal{K}}
\newcommand{\dd}{\,d}

\usepackage[colorlinks=true,linkcolor=blue,citecolor=darkgray,urlcolor=blue,bookmarks=true]{hyperref}
\hypersetup{
  pdftitle={A short proof of a Landen product identity for the Groetzsch ring function},
  pdfauthor={Junyi Hu}
}

\title{A short proof of a Landen product identity\\
for the Gr\"otzsch ring function}
\author[1]{Junyi Hu}
\affil[1]{Department of Automation, Tsinghua University}
\date{}

\begin{document}
\maketitle

\begin{abstract}
For $r\in(0,1)$ let $r_{0}=r$, $r_{n+1}=2\sqrt{r_{n}}/(1+r_{n})$ be the
ascending Landen sequence, and let
$\mu(r)=\frac{\pi}{2}\KK(r')/\KK(r)$ be the Gr\"otzsch ring function, where
$r'=\sqrt{1-r^{2}}$ and $\KK$ is the complete elliptic integral of the
first kind. We give a new, self-contained proof of the product identity
\begin{equation*}
\prod_{n=0}^{\infty}\bigl(1+r_{n}\bigr)^{2^{-n}}=r'\,e^{\mu(r')},
\end{equation*}
which is the equality case $a=1/2$ of a product theorem of Qiu and
Vuorinen. The proof uses the arithmetic--geometric mean, the theta duplication
identities, and a dyadic weighted $q$-series identity that is proved
directly by comparing coefficients. The symmetric closed form
$\log r+\mu(r)+\log P(r)=\log(rr')+\mu(r)+\mu(r')$ follows, and the
endpoint values $P(0^{+})=1$, $P(1^{-})=4$ are explained.
\end{abstract}

\section{Introduction and statement}\label{sec:intro}

For $r\in(0,1)$ write $r'=\sqrt{1-r^{2}}$ and let
\begin{equation}\label{eq:KK}
\KK(r)=\int_{0}^{\pi/2}\frac{\dd t}{\sqrt{1-r^{2}\sin^{2}t}}
\end{equation}
be the complete elliptic integral of the first kind. The decreasing
homeomorphism $\mu:(0,1)\to(0,\infty)$,
\begin{equation}\label{eq:mu}
\mu(r)=\frac{\pi}{2}\,\frac{\KK(r')}{\KK(r)},
\end{equation}
is the conformal modulus of the plane Gr\"otzsch ring
$\mathbb{B}^{2}\setminus[0,r]$ \cite{avv1997conformal}. The
\emph{ascending Landen sequence} of $r$ is defined by $r_{0}=r$ and
\begin{equation}\label{eq:landen}
r_{n+1}=\frac{2\sqrt{r_{n}}}{1+r_{n}},\qquad n=0,1,2,\dots;
\end{equation}
it satisfies $\mu(r_{n})=2^{-n}\mu(r)$ and converges to $1$
quadratically \cite{avv1997conformal,borwein1987pi}. The infinite
product
\begin{equation}\label{eq:P}
P(r)=\prod_{n=0}^{\infty}\bigl(1+r_{n}\bigr)^{2^{-n}}
\end{equation}
converges rapidly, since $1-r_{n}=O\bigl((1-r_{n-1})^{2}/8\bigr)$.

Qiu and Vuorinen proved the following product theorem
\cite[Theorem~1.10]{qiu1999infinite} (see also the survey
\cite[3.32]{vuorinen1997geometric}): for $a\in(0,1/2]$, with the
generalized Gr\"otzsch function $\mu_{a}$ and the Ramanujan constant
$R(a)$,
\begin{equation}\label{eq:qv}
P(r')\ \leq\ \exp\bigl(\mu_{a}(r)+\log r\bigr)\ \leq\
\frac{e^{R(a)}}{16}\,P(r'),
\end{equation}
with equality in both places if and only if $a=1/2$. Since
$\mu_{1/2}=\mu$ and $R(1/2)=\log 16$, the equality case reads
\begin{equation}\label{eq:main}
P(r)=r'\,e^{\mu(r')},\qquad r\in(0,1).
\end{equation}
The quantity $r'e^{\mu(r')}$ itself is classical: Jacobi's 1829 formula
\cite{jacobi1829fundamenta} (cf.\ \cite[p.~91]{avv1997conformal})
expresses $e^{\mu(r)+\log r}$ as the $q$-product
$4\prod_{n\geq1}\bigl((1+q^{2n})/(1+q^{2n-1})\bigr)^{4}$ with
$q=e^{-2\mu(r)}$.

The purpose of this note is to give a new, short and self-contained
proof of \eqref{eq:main}. The original proof of \eqref{eq:qv} in
\cite{qiu1999infinite} is based on infinite products and inequalities
for normalized quotients of hypergeometric functions. Our proof is
different: we telescope the product against the arithmetic--geometric
mean (AGM) iteration realized by the theta duplication identities, and
then settle the resulting identity by an elementary coefficient
comparison. The same computation yields a symmetric closed form for
the quantity in \eqref{eq:qv}:

\begin{theorem}\label{thm:main}
For $r\in(0,1)$ the product \eqref{eq:P} satisfies \eqref{eq:main}.
Consequently the function $G(r):=\log r+\mu(r)+\log P(r)$ is symmetric
in $r$ and $r'$:
\begin{equation}\label{eq:G}
G(r)=\log\bigl(rr'\bigr)+\mu(r)+\mu(r')=G(r').
\end{equation}
\end{theorem}

The key auxiliary identity, which may be of independent interest, is
the following dyadic weighted evaluation of $\theta_{3}$ at the points
$2^{m}\tau$.

\begin{lemma}[Dyadic theta identity]\label{lem:dyadic}
Let $\theta_{3}(\tau)=\sum_{k\in\mathbb{Z}}q^{k^{2}}$,
$q=e^{i\pi\tau}$, $\operatorname{Im}\tau>0$. Then
\begin{equation}\label{eq:dyadic}
\sum_{m=1}^{\infty}2^{1-m}\log\theta_{3}(2^{m}\tau)
=2\sum_{k=1}^{\infty}\log\bigl(1-q^{2k}\bigr)
+4\sum_{k=1}^{\infty}\log\bigl(1+q^{2k}\bigr),
\end{equation}
where on the left the real logarithm is taken for
$\tau\in i\mathbb{R}^{+}$ and the series converges absolutely.
\end{lemma}

\section{Proofs}\label{sec:proof}

\subsection{Theta parametrization}

We use the classical theta constants with none
$q_{\varsigma}=e^{i\pi\varsigma}$, $\operatorname{Im}\varsigma>0$, and
the modular lambda function
$\lambda(\varsigma)=\theta_{2}(\varsigma)^{4}/\theta_{3}(\varsigma)^{4}$.
Fix $r\in(0,1)$. For $\tau=it$ with $t=2\mu(r)/\pi$ one has
$\sqrt{\lambda(\tau)}=\theta_{2}(\tau)^{2}/\theta_{3}(\tau)^{2}=r$.
Set $\sigma=-1/\tau=i/t$. Since $\lambda(-1/\tau)=1-\lambda(\tau)$
\cite[21.51]{whittaker1927course},
\begin{equation}\label{eq:xsigma}
r=\frac{\theta_{4}(\sigma)^{2}}{\theta_{3}(\sigma)^{2}},\qquad
r'=\frac{\theta_{2}(\sigma)^{2}}{\theta_{3}(\sigma)^{2}},\qquad
q:=e^{i\pi\sigma}=e^{-\pi/t}\in(0,1),
\end{equation}
all theta values being positive on the imaginary axis, and
\begin{equation}\label{eq:murp}
\mu(r')=\frac{\pi\KK(r)}{2\KK(r')}=\frac{\pi}{2t}=-\frac{1}{2}\log q .
\end{equation}

\subsection{Duplication identities and the AGM}

For $m\geq0$ define
\begin{equation}\label{eq:agmab}
a_{m}=\frac{\theta_{3}(2^{m}\sigma)^{2}}{\theta_{3}(\sigma)^{2}},
\qquad
b_{m}=\frac{\theta_{4}(2^{m}\sigma)^{2}}{\theta_{3}(\sigma)^{2}} .
\end{equation}
The classical duplication identities
\cite[21.51]{whittaker1927course}
\begin{equation}\label{eq:duplication}
2\,\theta_{3}(2\varsigma)^{2}=\theta_{3}(\varsigma)^{2}+\theta_{4}(\varsigma)^{2},
\qquad
\theta_{4}(2\varsigma)^{2}=\theta_{3}(\varsigma)\,\theta_{4}(\varsigma),
\end{equation}
applied with $\varsigma=2^{m}\sigma$, give
$a_{m+1}=(a_{m}+b_{m})/2$ and $b_{m+1}=\sqrt{a_{m}b_{m}}$; since
$a_{0}=1$ and $b_{0}=r$, the pair $(a_{m},b_{m})$ is the
arithmetic--geometric iteration of $(1,r)$, and $a_{m}$ decreases to
Gauss's limit $M(1,r)=\pi/(2\KK(r'))>0$
\cite[Ch.~1]{borwein1987pi}. Put
$x_{m}:=b_{m}/a_{m}=\theta_{4}(2^{m}\sigma)^{2}/\theta_{3}(2^{m}\sigma)^{2}$.
Then $x_{0}=r$ and, by \eqref{eq:duplication},
\begin{equation*}
\frac{2\sqrt{x_{m}}}{1+x_{m}}
=\frac{2\,\theta_{3}(2^{m}\sigma)\,\theta_{4}(2^{m}\sigma)}
{\theta_{3}(2^{m}\sigma)^{2}+\theta_{4}(2^{m}\sigma)^{2}}
=\frac{\theta_{4}(2^{m+1}\sigma)^{2}}{\theta_{3}(2^{m+1}\sigma)^{2}}
=x_{m+1},
\end{equation*}
so $(x_{m})$ is precisely the Landen sequence \eqref{eq:landen};
moreover
\begin{equation}\label{eq:termwise}
1+x_{m}=\frac{a_{m}+b_{m}}{a_{m}}=\frac{2a_{m+1}}{a_{m}} .
\end{equation}

\subsection{Telescoping}

Summing the elementary identity
\begin{equation*}
\sum_{n=0}^{N}2^{-n}\bigl(\log a_{n+1}-\log a_{n}\bigr)
=\sum_{m=1}^{N}2^{-m}\log a_{m}+2^{-N}\log a_{N+1}
\qquad(a_{0}=1)
\end{equation*}
and letting $N\to\infty$, so that $2^{-N}\log a_{N+1}\to0$, we obtain
from \eqref{eq:termwise}
\begin{equation}\label{eq:telescope}
\log P(r)=\sum_{n=0}^{\infty}2^{-n}\log\frac{2a_{n+1}}{a_{n}}
=2\log2+\sum_{m=1}^{\infty}2^{1-m}\log\theta_{3}(2^{m}\sigma)
-2\log\theta_{3}(\sigma);
\end{equation}
the last series converges since
$\theta_{3}(2^{m}\sigma)=1+O(q^{2^{m}})$.

\subsection{Proof of the dyadic theta identity}

By Jacobi's triple product,
$\theta_{3}(\varsigma)=\prod_{k\geq1}(1-q_{\varsigma}^{2k})
(1+q_{\varsigma}^{2k-1})^{2}$ with $q_{\varsigma}=e^{i\pi\varsigma}$
\cite[21.3]{whittaker1927course}; expanding every logarithm into its
Taylor series turns both sides of \eqref{eq:dyadic} into absolutely
convergent double series in $q$ (the coefficients grow at most
polynomially in the exponent), so it suffices to compare the
coefficient $[q^{e}]$ for each $e\geq1$. Write $e=2^{v}o$ with $o$
odd and abbreviate $\sigma_{-1}(o)=\sum_{d\mid o}d^{-1}$.

On the right-hand side of \eqref{eq:dyadic} only even exponents
occur; the exponent has the form $e=2kj$, and writing $j=2^{s}d$ with
$d\mid o$, the constraint $k\geq1$ forces $0\leq s\leq v-1$. The
factor contributed by $2\log(1-y)+4\log(1+y)$ at $y^{j}$ is
$\bigl(-2+4(-1)^{j+1}\bigr)/j$, which equals $2/j$ for $s=0$ and
$-6/j$ for $s\geq1$. Hence, for $v\geq1$,
\begin{equation*}
[q^{e}]\,(\mathrm{RHS})
=\sigma_{-1}(o)\biggl(2-6\sum_{s=1}^{v-1}2^{-s}\biggr)
=\sigma_{-1}(o)\bigl(-4+3\cdot2^{2-v}\bigr).
\end{equation*}

On the left-hand side, the expansion of $\log(1-q^{2^{m+1}k})$
contributes, for each $1\leq m\leq v-1$, the amount
\begin{equation*}
-2^{1-m}\sigma_{-1}\bigl(2^{v-m-1}\bigr)\sigma_{-1}(o)
=-2^{1-m}\bigl(2-2^{m+1-v}\bigr)\sigma_{-1}(o),
\end{equation*}
which sums over $m$ to
$-\sigma_{-1}(o)\bigl(4-2^{3-v}-(v-1)2^{2-v}\bigr)$. The expansion of
$2\log(1+q^{2^{m}(2k-1)})$ contributes for each $1\leq m\leq v$: here
the cofactor $j$ must absorb the full dyadic part $2^{v-m}$ of
$e/2^{m}$ (so that $e/2^{m}j$ is odd), giving the amount
$\sigma_{-1}(o)\,2^{2-v}(2-v)$ in total. The sum is again
$\sigma_{-1}(o)\bigl(-4+3\cdot2^{2-v}\bigr)$, and for odd $e$ both
sides vanish. Lemma~\ref{lem:dyadic} is proved. \qed

\subsection{Assembly}

Inserting \eqref{eq:dyadic} into \eqref{eq:telescope} and using
Jacobi's product for $\theta_{2}$,
$\theta_{2}(\sigma)=2q^{1/4}\prod_{k\geq1}(1-q^{2k})(1+q^{2k})^{2}$
\cite[21.3]{whittaker1927course}, gives
\begin{align*}
\log P(r)
&=2\log2+2\sum_{k\geq1}\log\bigl(1-q^{2k}\bigr)
+4\sum_{k\geq1}\log\bigl(1+q^{2k}\bigr)-2\log\theta_{3}(\sigma)\\
&=2\log\theta_{2}(\sigma)-\frac{1}{2}\log q-2\log\theta_{3}(\sigma).
\end{align*}
By \eqref{eq:xsigma} and \eqref{eq:murp} the right-hand side equals
$\log r'+\mu(r')$, which is \eqref{eq:main}. The symmetry
\eqref{eq:G} follows at once, and Theorem~\ref{thm:main} is proved.
\qed

\section{Remarks}\label{sec:remarks}

\begin{remark}[Endpoint values]\label{rem:endpoints}
The closed form \eqref{eq:main} explains the endpoint values of the
product: $r'e^{\mu(r')}\to1$ as $r\to0^{+}$, while
$\mu(r)=\log(4/r)+O(r^{2}\log r)$ as $r\to0^{+}$
\cite[p.~82]{avv1997conformal} gives $r'e^{\mu(r')}\to4$ as
$r\to1^{-}$; that is, $P(0^{+})=1$ and $P(1^{-})=4$.
\end{remark}

\begin{remark}[Numerical check]\label{rem:numeric}
A numerical sample of \eqref{eq:G}:
\begin{equation*}
G\!\left(\tfrac12\right)
=\log\frac{\sqrt3}{4}+\mu\!\left(\tfrac12\right)
+\mu\!\left(\tfrac{\sqrt3}{2}\right)
=2.4003641589570148359\cdots,
\end{equation*}
in agreement with a direct evaluation of the product \eqref{eq:P}
(truncated after $n=150$) to more than $40$ digits.
\end{remark}

\begin{remark}[Relation to Jacobi's $q$-product]\label{rem:jacobi}
The proof above identifies the Landen product $P(r)$ with Jacobi's
$q$-product directly: the dyadic identity \eqref{eq:dyadic} is
precisely the statement that the weighted $\theta_{3}$-sum built from
the Landen sequence collapses to the triple-product expansion of
$\theta_{2}$. In this sense \eqref{eq:main} is the Landen-product
form of Jacobi's 1829 evaluation of $e^{\mu(r)+\log r}$
\cite{jacobi1829fundamenta,avv1997conformal}.
\end{remark}

\begin{remark}[On the case $a<1/2$]\label{rem:open}
It is natural to ask whether the coefficient comparison in
Lemma~\ref{lem:dyadic} has a hypergeometric $q$-analogue explaining
\eqref{eq:qv} for $a<1/2$. The answer appears to be twofold: the
hypergeometric coefficient comparison does exist---but is already in
the literature---while a literal $q$-analogue of \eqref{eq:dyadic}
apparently does not. For the first point, write
$D_{a}(r):=\mu_{a}(r)-\mu(r)$. By \eqref{eq:main} the chain
\eqref{eq:qv} is equivalent to $0\le D_{a}(r)\le R(a)-\log16$, and
the standard derivative formula
\begin{equation*}
\mu_{a}'(r)=-\frac{1}{r(1-r^{2})\,F(a,1-a;1;r^{2})^{2}}
\end{equation*}
($F$ the Gaussian hypergeometric function) shows that the sign of
$D_{a}'$ is decided by comparing the coefficients
\begin{equation*}
c_{n}(a)=\frac{(a)_{n}(1-a)_{n}}{(n!)^{2}}
=\frac{1}{(n!)^{2}}\prod_{k=0}^{n-1}\bigl(k(k+1)+a(1-a)\bigr)
\end{equation*}
of $F(a,1-a;1;x)$, which increase with $a$ on $(0,1/2]$
\cite[Theorem~7.2(3)]{aqvv2000generalized}. Hence $D_{a}$ is strictly
decreasing; since $D_{a}(0^{+})=(R(a)-\log16)/2$ and $D_{a}(1^{-})=0$
\cite[Theorem~5.5(2)]{aqvv2000generalized}, this gives the sharp form
of the upper bound in \eqref{eq:qv},
\begin{equation*}
\exp\bigl(\mu_{a}(r)+\log r\bigr)\ \le\
\frac{e^{R(a)/2}}{4}\,P(r'),\qquad a\in(0,1/2],
\end{equation*}
with the best constant $e^{R(a)/2}/4$, the form in which the
Qiu--Vuorinen theorem is stated in \cite[(1.5)]{bao2022notes}. The
same coefficient method is carried much further in
\cite{bao2022notes}, where $\mu_{a}(r)+\log
r=R(a)/2-\sum_{n\ge1}\theta_{n}r^{2n}$ with computable positive
coefficients $\theta_{n}$.

For the second point, there is numerical evidence against any triadic
analogue: with the cubic theta function
$a(q)=\sum_{j,k\in\mathbb{Z}}q^{j^{2}+jk+k^{2}}$, the formal series
$\sum_{m\ge1}3^{1-m}\log a(q^{3^{m}})$ does not collapse to a
Lambert-type series (its coefficients grow with alternating sign).
The dyadic identity \eqref{eq:dyadic} thus seems to be a
signature-$2$ phenomenon, relying both on Jacobi's triple product for
$\theta_{3}$ and on the exact match between the weights $2^{1-m}$ and
the $2$-adic valuation of the exponents. This fits the shape of
\eqref{eq:qv}: for $a<1/2$ both inequalities are strict, so there is
no closed form awaiting an identity; equality holds only at $a=1/2$,
which is precisely \eqref{eq:main}.
\end{remark}

\paragraph{Acknowledgements.}
The identity \eqref{eq:main} was rediscovered independently by the
authors; the priority of Qiu and Vuorinen \cite{qiu1999infinite} is
acknowledged.


\begin{thebibliography}{9}

\bibitem{aqvv2000generalized}
G.~D. Anderson, S.-L. Qiu, M.~K. Vamanamurthy, and M. Vuorinen,
Generalized elliptic integrals and modular equations,
\emph{Pacific J. Math.} \textbf{192} (2000), no.~1, 1--37.

\bibitem{avv1997conformal}
G.~D. Anderson, M.~K. Vamanamurthy, and M.~Vuorinen,
\emph{Conformal Invariants, Inequalities, and Quasiconformal Maps},
John Wiley \& Sons, New York, 1997.

\bibitem{bao2022notes}
Q. Bao and M.-K. Wang,
Notes on generalized Gr\"otzsch ring function and generalized
Hersch--Pfluger distortion function, preprint, 2022,
arXiv:2202.09758.

\bibitem{borwein1987pi}
J.~M. Borwein and P.~B. Borwein,
\emph{Pi and the AGM: A Study in Analytic Number Theory and
Computational Complexity}, John Wiley \& Sons, New York, 1987.

\bibitem{jacobi1829fundamenta}
C.~G.~J. Jacobi,
\emph{Fundamenta Nova Theoriae Functionum Ellipticarum},
K\"onigsberg, 1829.

\bibitem{qiu1999infinite}
S.-L. Qiu and M. Vuorinen,
Infinite products and normalized quotients of hypergeometric
functions, \emph{SIAM J. Math. Anal.} \textbf{30} (1999), no.~5,
1057--1075.

\bibitem{vuorinen1997geometric}
M. Vuorinen,
Geometric properties of quasiconformal maps and special functions,
I--III, \emph{Bull. Soc. Sci. Lett. \L\'od\'z S\'er. Rech. D\'eform.}
\textbf{24} (1997), 7--58; also arXiv:math/0703687.

\bibitem{whittaker1927course}
E.~T. Whittaker and G.~N. Watson,
\emph{A Course of Modern Analysis}, 4th ed., Cambridge University
Press, Cambridge, 1927.

\end{thebibliography}
\end{document}